\documentclass[lean]{AFM}

\usepackage{amsmath, amsfonts, amsthm, framed}
\usepackage{xspace}
\usepackage{newunicodechar}
\usepackage{enumitem}
\usepackage{fontawesome}
\usepackage{soul}
\usepackage{cancel}
\usepackage{xspace}
\usepackage{eucal}
\usepackage{colonequals}
\definecolor{darkergreen}{rgb}{0.0, 0.5, 0.0}
\definecolor{Blue}{RGB}{0,0,148}
\definecolor{keywordcolor}{rgb}{0.7, 0.1, 0.1}   % red
\definecolor{commentcolor}{rgb}{0.4, 0.4, 0.4}   % grey
\definecolor{symbolcolor}{rgb}{0, 0, 0.8}    % blue
\definecolor{tacticcolor}{rgb}{0, 0, 0.8}    % blue
\definecolor{sortcolor}{rgb}{0.1, 0.5, 0.1}      % green
\theoremstyle{plain}
\theoremstyle{definition}
\newcommand{\CC}{\mathbb{C}}
\newcommand{\RR}{\mathbb{R}}
\newcommand{\ZZ}{\mathbb{Z}}
\newcommand{\NN}{\mathbb{N}}

\renewcommand{\le}{\leq}
\renewcommand{\leq}{\leqslant}
\renewcommand{\ge}{\geq}
\renewcommand{\ge}{\geqslant}

\newcommand{\docs}[1]{%
\href{https://leanprover-community.github.io/mathlib4_docs/find/?pattern=#1\#doc}
{\texttt{\detokenize{#1}}}}

\title[Formalizing zeta and L-functions in Lean]{Formalizing zeta and L-functions in Lean}

\author[D. Loeffler and M. Stoll]{David Loeffler and Michael Stoll}

\authorinfo[D. Loeffler]{UniDistance Suisse}{david.loeffler@unidistance.ch} 
\authorinfo[M. Stoll]{Mathematisches Institut, Universität Bayreuth}{michael.stoll@uni-bayreuth.de} 
\VOLUME{1}
\YEAR{2025}
\firstpage{43}
\DOI{https://doi.org/10.46298/afm.15328}
\receiveddate{March 6, 2025}
\finaldate{June 13, 2025}
\accepteddate{June 23, 2025} % update this



\begin{abstract}
The Riemann zeta function, and more generally the $L$-functions of Dirichlet characters, are among the central objects of study in analytic number theory. We report on a project to formalize the theory of these objects in Lean's `Mathlib' library, including a proof of Dirichlet's theorem on primes in arithmetic progressions and a formal statement of the Riemann hypothesis.
\end{abstract}

\support{D.L. gratefully acknowledges the support of the European Research Council through the Horizon 2020 Excellent Science programme (Consolidator Grant ``ShimBSD: Shimura varieties and the BSD conjecture'', grant ID 101001051). M.S. thanks the Institute for Mathematical Research (FIM) at ETH Zürich for enabling him to spend a week at ETH.}

\begin{document}

\section{Introduction}
 \subsection{Mathematical background}

  The \emph{Riemann zeta function} is the unique holomorphic function $\zeta : \CC - \{1\} \to \CC$ whose restriction to $\{ s : \Re(s) > 1\}$ agrees with the sum of the series $\sum_{n = 1}^{\infty} \tfrac{1}{n^s}$.

  This function has a long history, going back at least to Euler, who proved two important results about this function: firstly, the remarkable formula $\zeta(2) = \tfrac{\pi^2}{6}$ (solving the so-called ``Basel problem''); and secondly, the product formula
  \[ \tag{\dag} \zeta(s) = \prod_{\text{$p$ prime}} \left(1- \tfrac{1}{p^s} \right)^{-1}, \]
  which points to a deep connection between $\zeta(s)$ and the prime numbers.

  The contribution of Riemann, in his epochal 1859 paper \cite{Ri59}, was to to apply methods of complex analysis to the zeta function -- firstly showing it has analytic continuation beyond $\{\Re(s) > 1\}$, and secondly showing that the location of the zeroes and poles of this extended function controls the distribution of the prime numbers. Using this function, Riemann sketched a path to proving a conjecture of Gauss, namely that as $X \to \infty$, the number of primes $\le X$ is asymptotically $\tfrac{X}{\log X}$. Riemann's program was successfully completed out by Hadamard and de la Vall\'ee Poussin in 1896, proving the asymptotic formula now known as the \emph{Prime Number Theorem}.

  Riemann's zeta-function can be seen as just one example of a wider class of functions: Dirichlet $L$-functions $L(\chi, s)$, attached to a Dirichlet character $\chi$ modulo $N$ for some $N \ge 1$ (i.e.~a group homomorphism $(\ZZ / N\ZZ)^\times \to \CC^\times$). Here $L(\chi, s)$ is defined as the analytic continuation of the function on $\{ \Re(s) > 1\}$ defined by the series
  \[ \sum_{\substack{n \ge 1 \\ \gcd(n, N) = 1}} \frac{\chi(n)}{n^s}. \]
  Dirichlet introduced these functions in order to study how the primes are distributed among the residue classes modulo $N$; they played a central role in his proof that for any $a$ with $\gcd(a, N) = 1$, there exist infinitely many primes $p$ with $p \equiv a \bmod N$.

\subsection{Results formalized in this project}

  The aim of this paper is to report on a project to formalize the definitions and key properties of zeta and $L$-functions, as well as a proof of Dirichlet's theorem mentioned above, in the Lean theorem prover, and to contribute these formalizations to Mathlib, Lean's mathematics library.

  The definitions and theorems about the zeta-function added to Mathlib in this project include\footnote{The name of each declaration is a link to its entry in the Mathlib reference manual. Observe that the naming follows Lean's convention that names of terms (such as functions, or proofs of theorems) start with small letters, and names of types (such as statements of conjectures) with capitals.}:
  \begin{itemize}
  \item \docs{riemannZeta}: a definition of the Riemann zeta function.

 % \item \docs{differentiableAt_riemannZeta}: proof that $\zeta$ is holomorphic away from $s = 1$.

  %\item \docs{zeta_eq_tsum_one_div_nat_cpow}: proof that it agrees with the sum of the series $\sum_{n  = 1}^\infty n^{-s}$ for $\Re(s) > 1$.
  \item \docs{riemannZeta_two}: the formula $\zeta(2) = \tfrac{\pi^2}{6}$.
  \item \docs{riemannZeta_eulerProduct}: Euler's product formula ($\dag$).
  \item \docs{riemannZeta_one_sub}: the functional equation relating $\zeta(s)$ to $\zeta(1 - s)$.

  %\item \docs{riemannZeta_two_mul_nat}: formula for $\zeta(2 k)$ for $k \ge 1 \in \mathbb{N}$ in terms of Bernoulli numbers (and the special case \docs{riemannZeta_two} for the formula $\zeta(2) = \tfrac{\pi^2}{6}$).
%  \item \docs{tendsto_riemannZeta_sub_one_div}: the formula
%  \[ \zeta(s) = \tfrac{1}{s - 1} + \gamma + o(1) \quad \text{as $s \to 1$},\]
%  where $\gamma$ is the Euler--Mascheroni constant $\lim_{X \to \infty} \sum_{n = 1}^X \tfrac{1}{n} - \log X$.

  \item \docs{riemannZeta_ne_zero_of_one_le_re}: proof that $\zeta(s) \ne 0$ for all $s$ in the closed half-plane $\Re(s) \ge 1$ (including on the boundary $\Re(s) = 1$).

  \item \docs{RiemannHypothesis}: a formal \emph{statement} of the ``Riemann hypothesis'', i.e.~that all zeroes of $\zeta(s)$, other than the trivial zeroes at negative even integers, have $\Re(s) = \tfrac{1}{2}$.
  \end{itemize}

  There are analogous results for Dirichlet characters, including the non-vanishing on  $\{ s : \Re(s) \ge 1\}$. From this we deduce \docs{Nat.setOf_prime_and_eq_mod_infinite}, Dirichlet's theorem on primes in arithmetic progressions.

 \subsection{Related work}

  \begin{itemize}
   \item There is an ``elementary'' proof of the Prime Number Theorem (avoiding complex analysis and zeta-functions completely) due to Erd\H{o}s and Selberg, and another shorter but less elementary analytic proof due to Newman. These have been formalized in several proof systems; see the introduction of \cite{PNT+} for references. However, these arguments do not give such good quantitative bounds on the error term as the classical analytic argument, and largely side-step developing the theory of zeta and $L$-functions, which is an interesting goal in its own right.

   \item Gomes and Kontorovich \cite{GK20} formalized in Lean a statement \emph{equivalent} to the Riemann hypothesis, formulated using the Dirichlet $\eta$ function
   \[ \eta(s) = 1 - \tfrac{1}{2^s} + \tfrac{1}{3^s} - \dots \]
   which is (conditionally) convergent for $\Re(s) > 0$, and sums to $(1 - 2^{1-s}) \zeta(s)$. It would be interesting to formalize the theorem that the version of the Riemann hypothesis formalized by Gomes--Kontorovich is equivalent to ours, but we have not attempted to do this.

   \item Manuel Eberl and his collaborators have formalized a substantial part of the theory of Riemann zeta and $L$-functions in the ``Isabelle'' theorem prover \cite{Eb19}. Building on this, Song and Yao \cite{SY24} recently announced an Isabelle formalization of the prime number theorem with the ``classical'' error term $O(x \exp(-C \sqrt{\log x}))$).

   While the main results of this Isabelle-based project are comparable to (and indeed go considerably further than) the results of the Lean formalization described here, there are major differences in methodology, particularly in the proof of the functional equation. For this result, Riemann gave two proofs in \cite{Ri59}. The Isabelle proof is based on contour integrals, while we follow the theta-function proof. As we shall see below, this involved the formalization of a range of interesting results along the way, and has many natural generalizations, such as in the theory of modular forms.

   \item Building on the formalizations described in this paper, the ``PrimeNumberTheorem+'' project led by Alex Kontorovich and Terry Tao \cite{PNT+} has recently formalized in Lean a proof of the Prime Number Theorem (via the Wiener--Ikehara theorem), which will hopefully be merged into Mathlib in the near future. Future goals of the PNT+ project include formalizing an explicit bound on the error term, and extending the results to counting primes in arithmetic progressions, strengthening Dirichlet's theorem by giving an asymptotic for the number of primes $\le X$ in a given congruence class.

  \end{itemize}

\section{Implementing $L$-series in Mathlib}

  The number-theoretic significance of the Riemann zeta function, Dirichlet $L$-functions
  and numerous other $L$-functions is related to the sequence of coefficients in their
  expansion as Dirichlet series and its properties (e.g., via an expression as an
  Euler product). So an implementation of Dirichlet series (which we will call $L$-series below)
  and their relevant properties is a prerequisite for applications like the Prime Number Theorem or Dirichlet's Theorem.

  At the beginning of this project, a rudimentary implementation of $L$-series in Mathlib
  existed, which associated a function $\CC \to \CC$ to an ``arithmetic function''~$f$;
  arithmetic functions are implemented as functions $f \colon \NN \to R$ with $f(0) = 0$
  (here $R$ is some semiring, which needs to have a coercion to~$\CC$ in our context).
  This made it possible to write \texttt{LSeries}~$\zeta$ or \texttt{LSeries}~$\mu$
  for the Dirichlet series associated to the arithmetic functions~$\zeta$ (taking the
  value~$1$ everywhere outside zero) or the Möbius function~$\mu$; but, for example,
  \texttt{LSeries}~$\chi$ with a Dirichlet character~$\chi$, or \texttt{LSeries}~$\Lambda$
  with the von Mangoldt function~$\Lambda$, did not type-check because there are no
  coercions from \texttt{DirichletCharacter}~$\CC$~$n$ or from \texttt{ArithmeticFunction}~$\RR$
  to \texttt{ArithmeticFunction}~$\CC$. 
  
  We found it was not easily possible to set
  up a coercion from \texttt{ArithmeticFunction}~$R$ to
  \texttt{ArithmeticFunction}~$R'$ whenever there is a coercion from $R$ to~$R'$, because the necessary coercion instance would have too many free parameters; 
  and hard-coding the relevant cases $\NN \hookrightarrow \ZZ \hookrightarrow \RR \hookrightarrow \CC$ led to a considerable amount of code duplication.
  See the discussion on the Lean Zulip
  chat\footnote{\url{https://leanprover.zulipchat.com/\#narrow/channel/144837-PR-reviews/topic/.2310725.20and.20.2310728.3A.20L-series/near/422590170}}
  for details on this and also on how we arrived at the current design, which is described below.
  
  To avoid these problems, we changed the implementation of $L$-series to not use
  the type of arithmetic functions. A first attempt using the argument type $\NN+ \to \CC$
  (where $\NN+$ is the subtype of~$\NN$ consisting of positive natural numbers)
  was abandoned in favour of the approach described below because it turned out
  that the necessity of going back and forth between $\NN$ and~$\NN+$ was causing
  some friction; also, some API for~$\NN+$ was missing.

  The design that we finally settled on is to use functions $f \colon \NN \to \CC$
  and simply ignore the value of~$f$ at zero, so that
  \[ \operatorname{\texttt{LSeries}} \; f \; s = \sum_{n \ne 0} \frac{f(n)}{n^s} \]
  whenever the series on the right converges absolutely (and the value is zero
  otherwise by the Mathlib convention for topological sums); see \docs{LSeries}. To solve the problem
  with \texttt{LSeries}~$\chi$ etc., we introduce the notation ${}^\nearrow f$,
  which coerces $f \colon A \to R$ to a function from~$\NN$ to~$\CC$ when
  coercions $\NN \to A$ and $R \to \CC$ are available. We introduce predicates
  \docs{LSeriesHasSum} and~\docs{LSeriesSummable} to be able to talk about
  when an $L$-series converges.

  It turns out that with this set-up, the formalization of operations with
  and properties of $L$-series proceeds quite smoothly. This includes, for example,
  the Euler product representation of Dirichlet $L$-series.

  The basic theory of $L$-series is developed in about~1400 lines of Lean code
  spread over various files under \href{https://github.com/leanprover-community/mathlib4/tree/master/Mathlib/NumberTheory/LSeries}{\texttt{Mathlib/NumberTheory/LSeries}}, and the derivation
  of the Euler product representation (in files under \href{https://github.com/leanprover-community/mathlib4/tree/master/Mathlib/NumberTheory/EulerProduct}{\texttt{Mathlib/NumberTheory/EulerProduct}})
  takes another 500~lines of code.

\section{From Fourier analysis to zeta functions}

 \subsection{Fourier analysis}

  The key analytic input for our treatment of zeta and $L$-functions is the theory of Fourier series, expressing functions $\RR / \ZZ \to \CC$ (i.e.~functions on $\RR$ which are periodic with period 1) in terms of the Fourier basis functions $x \mapsto e^{2 \pi i n x}$ for $n \in \ZZ$.

  At the outset of this project, Mathlib already contained some substantial results in Fourier theory (mostly contributed by Heather Macbeth), leading up to a proof that the Fourier basis functions are an orthonormal basis of the Hilbert space $L^2(\RR / \ZZ)$ of square-integrable functions\footnote{More precisely, equivalence classes of functions differing on a set of zero measure.} on $\RR/\ZZ$ (\docs{hasSum_fourier_series_L2}). For our applications, we needed to formalize several more additional results in Fourier theory:

  \begin{itemize}
  \item \emph{Uniform convergence of Fourier series}: if $f$ is a continuous function on $\RR / \ZZ$ such that $\sum_{n \in \ZZ} |c_n(f)| < \infty$, where $c_n(f)$ is the $n$-th Fourier coefficient $\int_0^1 e^{-2 \pi i n z} f(z)\, \mathrm{d}z$, then the Fourier series of $f$ converges absolutely and uniformly to $f$.

  \item \emph{Fourier transforms} for functions on the whole real line: definition and basic properties.

  \item \emph{Poisson's summation formula}: if $f : \RR \to \CC$ is continuous and absolutely integrable, and both $f$ and its Fourier transform $\hat f$ have sufficiently rapid decay at $\pm \infty$, then
  \[ \sum_{n \in \ZZ} f(n) = \sum_{n \in \ZZ} \hat{f}(n). \]
  \end{itemize}

  The uniform-convergence result for continuous functions follows from the $L^2$ version by arguing that if $\sum_{n \in \ZZ} |c_n(f)| < \infty$, then the Fourier series must converge uniformly to \emph{something}; and since uniform convergence implies $L^2$ convergence, and we know that the series converges in $L^2$ to $f$, it must in fact converge uniformly to $f$. This is formalized as \docs{hasSum_fourier_series_of_summable}.

  To deduce the Poisson summation formula, we show that for any sufficiently well-behaved function $f : \RR \to \CC$, the function
  \[ F(x) = \sum_{k \in \ZZ} f(x + k) \]
  is a continuous periodic function, and its Fourier coefficients are given by $c_n(F) = \hat{f}(n)$. So Poisson's summation formula amounts to the pointwise convergence of the Fourier series of $F$ at $0$, which is an application of the preceding theorem. This is formalized as \docs{Real.tsum_eq_tsum_fourierIntegral_of_rpow_decay}, assuming that $f$ is continuous with both $|f|$ and $|\hat{f}|$ decaying as $O(|x|^{-b})$ at infinity for some $b > 1$.

  The Fourier-analytic parts of this project amounted to approximately 1400 lines of code in Mathlib.

  \subsection*{Remark} Note that in this formulation we have decay assumptions on both $f$ and $\hat{f}$. This is easy to verify in our application, where both $f$ and $\hat{f}$ are Gaussian functions. However, in other applications it is advantageous to formulate all the hypotheses in terms of $f$ alone. Building on the contributions described above, S.~Gou\"ezel has recently added to Mathlib a proof that the Fourier transform of a Schwartz function is also a Schwartz function, from which it follows that the Poisson summation formula holds for any Schwartz function (\docs{SchwartzMap.tsum_eq_tsum_fourierIntegral}).

 \subsection{Jacobi's theta function}

  The \emph{Jacobi theta function} $\theta(\tau)$ is defined, for $\tau \in \CC$ with $\operatorname{Im}(\tau) > 0$, by the series
  \[ \theta(\tau) = \sum_{n \in \ZZ} e^{\pi i n^2 \tau}. \]
  In Mathlib this is defined as \docs{jacobiTheta}. Many properties of this function, such as holomorphy in $\tau$ and periodicity under $\tau \mapsto \tau + 2$, are straightforward from the definition. A much deeper symmetry relation is the transformation law (\docs{jacobiTheta_S_smul})
  \[
    \theta\left(\frac{-1}{\tau}\right) = \sqrt{-i\tau} \theta(\tau).
  \]
  This is derived from Poisson's summation formula applied to the Gaussian function $x \mapsto e^{\pi i x^2 \tau}$, whose Fourier transform is another Gaussian with $\tau$ replaced by $\tfrac{-1}{\tau}$; this is a version of the famous Gaussian integral $\int_{-\infty}^\infty e^{-x^2}\, \mathrm{d}x = \sqrt{\pi}$, recently contributed to Mathlib by S.~Gou\"ezel (\docs{integral_gaussian}). Since the Mellin transform $\int_0^\infty t^{s - 1} \frac{\theta(it) - 1}{2}\, \mathrm{d}t$ is $\pi^{-s} \Gamma(s) \zeta(2s)$, this gives the meromorphic continuation and functional equation of the zeta function.

  In order to extend this to Dirichlet character $L$-series, it is necessary to consider the more general two-variable theta function (formalized as \href{https://leanprover-community.github.io/mathlib4_docs/Mathlib/NumberTheory/ModularForms/JacobiTheta/TwoVariable.html}
  {\texttt{jacobiTheta$_\mathtt{2}$}})
  \[ \theta(z, \tau) = \sum_{n \in \ZZ} e^{2\pi i n z + \pi i n^2 \tau}.\]
  This satisfies a functional equation generalizing that of the one-variable function (and proved by the same method); and its Mellin transform gives the analytic continuation and functional equation of the even Hurwitz zeta function (\docs{HurwitzZeta.hurwitzZetaEven})
  \[
   \zeta_{\alpha}^{\text{ev}}(s) = \tfrac{1}{2}\sum_{\substack{n \in \ZZ \\ n + \alpha \ne 0}} \frac{1}{|n + \alpha|^s},\qquad \alpha \in \RR.
  \]
  By taking a suitable weighted sum over the values $\alpha = \tfrac{k}{N}$ for $k \in \ZZ / N\ZZ$, we can obtain the analytic continuation and functional equation for any \emph{even} Dirichlet character (satsifying $\chi(-1) = 1$).

  To include the case of \emph{odd} Dirichlet characters, we need to study yet another variant of the theta-series (\href{https://leanprover-community.github.io/mathlib4_docs/Mathlib/NumberTheory/ModularForms/JacobiTheta/TwoVariable.html}
    {\texttt{jacobiTheta$_\mathtt{2}'$}}), defined by
  \[
   \theta'(z, \tau) = \tfrac{\mathrm{d}}{\mathrm{d}z} \theta(z, \tau) = \sum_{n \in \ZZ} 2\pi i n\, e^{(2\pi i n z + \pi i n^2 \tau)}.
  \]
  This is related to the odd Hurwitz zeta function (\docs{HurwitzZeta.hurwitzZetaOdd}), defined by
  \[ \zeta_{\alpha}^{\text{odd}}(s) = \tfrac{1}{2}\sum_{\substack{n \in \ZZ \\ n + \alpha \ne 0}} \frac{\operatorname{sign}(n + \alpha)}{|n + \alpha|^s}, \]
  and as before we can obtain $L$-functions of odd Dirichlet characters by summing over $\alpha = \tfrac{k}{n}$.

  In the Mathlib implementation of the above results, the definition and properties of the Jacobi theta function is approximately 900 code lines, and the proof of the analytic continuation and functional equation of Dirichlet $L$-functions another 3300 lines. (This surprisingly large total is in part because several intermediate steps, e.g.~justifying the interchange of integration and summation required to relate $\int_0^\infty t^{s - 1} \frac{\theta(it) - 1}{2}\, \mathrm{d}t$ to $\sum_n n^{-s}$ when the latter converges, were done in considerably greater generality than absolutely required for this project, in order to enable later generalizations; see below.)

  \subsubsection*{Remarks} Since the above results on Poisson summation and Jacobi theta functions were added to Mathlib, they have found a surprising and unexpected real-world application: the transformation formula for theta functions was used in implementing verified algorithms for privacy-preserving computation for Amazon Web Services \cite{dM+24}.

 \subsection{Generalizations}

  The Jacobi theta function is an example of a modular form of weight $\tfrac{1}{2}$ (and level 4). The arguments used in our project for proving analytic continuation and functional equations for Riemann zeta and Dirichlet $L$-functions are special cases of more general arguments, which can be used to prove analogous results for $L$-functions of any modular form.

  Hence we have set up these arguments in an axiomatic framework which is intended to apply to any modular form $L$-function, following the arguments in \cite[Chapter 5]{DS05}. We introduce a concept we call a ``functional equation pair'', or ``FE-pair'' for short: this is a pair of locally integrable functions $f, g$ on $\RR_{\ge 0}$ such that
  \begin{itemize}
   \item $f(x)$ and $g(x)$ each have the form (constant) + (rapidly decaying term) as $t \to \infty$,
   \item $f (1 / x) = \varepsilon x^k g (x)$, for some positive real $k$ and complex $\varepsilon$.
  \end{itemize}
  The above definition is formalized as \docs{WeakFEPair} (with \docs{StrongFEPair} being the special case where the constant terms of $f, g$ at $\infty$ are both zero). We show that for any FE-pair, the Mellin transforms of $f$ and $g$ have meromorphic continuation (with poles and residues determined by the asymptotics of $f$ and $g$) and satisfy a functional equation relating values at $s$ and $k - s$ (\docs{WeakFEPair.functional_equation}).

  The basic example of an FE-pair is $f(x) = g(x) = \theta(ix)$ (with $k = \tfrac{1}{2}$); and the two-variable theta functions also give rise to FE-pairs, again with $k = \tfrac{1}{2}$. The formalization of modular form theory in Mathlib is ongoing; and we hope to use the FE-pair framework to formalize properties of $L$-functions of higher weight modular forms in future projects.

\section{Special values of $L$-functions}

 Lastly, we describe the formalization of Euler's ``Basel problem'' formula $\zeta(2) = \tfrac{\pi^2}{6}$, and its generalizations. The proof we formalized for these also uses Fourier analysis, although in a rather different fashion from above. It arises by studying the ``periodized Bernoulli function'' \docs{periodizedBernoulli}, i.e.~the unique function on $\RR / \ZZ$ whose restriction to $[0, 1)$ is the $k$-th Bernoulli polynomial $B_k(x)$. A standard computation shows that the $n$-th Fourier coefficient of this function is $\tfrac{-k!}{(2 \pi i n)^k}$. For even $k \ge 2$, considering the pointwise sum at $0$, and applying the results on Fourier-series convergence described above, gives the standard formula for $\zeta(2k)$ in terms of the $2k$-th Bernoulli number $B_{2k}(0)$ (\docs{riemannZeta_two_mul_nat}). This includes the Basel problem as a special case.

 This approach generalizes immediately to computing values of Hurwitz zeta functions. The result is most tidily stated in terms of values at \emph{negative} integers:
 \[ \zeta_{\alpha}(-k) = \frac{-B_{k+1}(\alpha)}{k+1} \quad \text{for $\alpha \in [0, 1]$ and $k \ge 1$}. \]
This is formalized as \docs{HurwitzZeta.hurwitzZeta_neg_nat}.

 A disadvantage of our approach is that we must exclude the case $k = 0$. In this case the above formula remains valid for $\alpha \in (0, 1)$; but our proof no longer applies, since the Fourier series concerned has $c_n(f) = -1/(2 \pi i n)$ and hence does not satisfy the condition $\sum_{n \in \ZZ} |c_n(f)| < \infty$. For instance, the Gregory--Leibniz formula $1 - \tfrac{1}{3} + \tfrac{1}{5} - \tfrac{1}{7} + \dots = \tfrac{\pi}{4}$, giving the value at $s = 1$ of the $L$-series of the unique nontrivial Dirichlet character modulo 4, is not a special case of our results\footnote{This particular formula is already in Mathlib as \docs{Real.tendsto_sum_pi_div_four}, but the method used does not seem to generalize to other Dirichlet characters. Since the original draft of this paper was written, X. Roblot has contributed to Mathlib a proof of Dedekind's class number formula for general number fields, \docs{NumberField.tendsto_sub_one_mul_dedekindZeta_nhdsGT}; applying this to quadratic fields, it should be relatively straightforward to deduce formulae for $L(\chi, 1)$ for any quadratic Dirichlet character.}. Including this case would require extending Mathlib's Fourier theory library by adding more general criteria for the Fourier series of a periodic function to converge to the function (such as the Dirichlet--Jordan criterion for functions of bounded variation); this would be an interesting project for the future.

\section{Pitfalls: totalizing functions}

 Mathematically, the natural domain of $\zeta(s)$ is the ``punctured complex plane'' $\CC - \{1\}$; there is no sensible definition of $\zeta(1)$ as an element of $\CC$. However, in formalizing proofs it is inconvenient to work with functions defined on subtypes; so the usual practice is to extend partially-defined functions to total functions by assigning \emph{junk values} at the bad points, such as the convention that $1 / 0 = 0$ used by many theorem provers (including Coq, Isabelle and Lean) \cite{div0}.

 Thus Mathlib's \docs{riemannZeta} function has a well-defined value for all complex numbers, including $1$. However, this junk value was not explicitly chosen but rather ``propagated'' from other choices of junk values at earlier states of the construction. This led to the following somewhat embarrassing episode. When the first author initially announced a formalization of the Riemann Hypothesis, the statement amounted to:

 \begin{quotation}
  ``If $s \in \CC$ is not a strictly negative even integer, and $\zeta(s) = 0$, then $\Re(s) = \tfrac{1}{2}$.''
 \end{quotation}

 As pointed out by K.~Buzzard and J.~Ellenberg (in discussions following a Twitter post by Buzzard about this work), there is an issue if $s = 1$. Our construction assigns some value to the expression \texttt{riemannZeta 1}; but it is far from obvious what this value is. If this junk value happened to be zero, then the conjecture we had formulated would not be the Riemann Hypothesis; it would be trivially false!

 Fortunately, we were able to show that the Mathlib definition of $\zeta(1)$ has a non-zero value: it is $\tfrac{1}{2}(\gamma - \log 4\pi)$, where $\gamma$ is the Euler--Mascheroni constant. This peculiar-looking formula arises because the zeta function is defined in three steps:
 \begin{enumerate}
 \item Firstly, we define the function
 \[ \Lambda_0(s) = \int_1^\infty (t^{s/2} + t^{(1-s)/2}) \frac{\theta(it) - 1}{2}\, \mathrm{d}t,\]
 which is holomorphic on all of $\CC$;
 \item we define $\Lambda(s) = \Lambda_0(s) - \tfrac{1}{s} - \tfrac{1}{1-s}$;
 \item finally we define $\zeta(s) = \frac{\pi^{s/2}}{\Gamma(s/2)} \Lambda(s)$ (up to a correction at $s = 0$).
 \end{enumerate}
 The ``junk value'' at $s = 1$ enters at step 2, since $\tfrac{1}{1-1}$ is defined to be 0. It follows that Mathlib's definition of $\zeta(1)$ is equal to the limit
 \[ \Lambda_0(1) - 1 = \lim_{s \to 1} \left(\zeta(s) - \frac{\pi^{s/2}}{(s - 1) \Gamma(s/2)}\right), \]
 and a lengthy and quite delicate computation, using the formula\footnote{This theorem (\docs{tendsto_riemannZeta_sub_one_div}) was in fact added to the project specifically for this purpose.} $\zeta(s) = \tfrac{1}{s - 1} + \gamma + o(1)$ and properties of the Gamma function, shows that this limit is $\tfrac{1}{2}(\gamma - \log 4\pi)$ as claimed. With this in hand, it is not difficult to formalize sufficiently good bounds on the constants involved to show that $\zeta(1) \ne 0$ (\docs{riemannZeta_one_ne_zero}); its actual value is about $-0.98$.

 \subsection*{Remark} Such oversights might have unfortunate consequences in the context of efforts to automate discovery of proofs via artificial intelligence. In this situation, mistakes in formalizing the \emph{statements} of conjectures might lead to mathematicians believing that a conjecture had been proved (or disproved), when in fact the system had proved or disproved a superficially similar but much easier statement, ``exploiting'' an oversight in formulating the conjecture in some mathematically uninteresting corner case.

\section{Applications to primes in arithmetic progressions}

  We can use the general results on $L$-series and the specific results on the
  analytic continuation of Dirichlet $L$-series to obtain a proof of Dirichlet's Theorem
  on primes in arithmetic progressions, following the classical analytic proof.

  One key ingredient that is needed here is a proof that the Dirichlet $L$-functions do not vanish on the \emph{closed} right half-plane $\Re(s) \ge 1$. When $s \neq 1$ or the Dirichlet character is not quadratic, this follows fairly easily from the Euler product representation and the non-negativity of a certain trigonometric polynomial. The remaining case, $L(\chi, 1) \neq 0$ when $\chi$ is a quadratic Dirichlet character, requires a more subtle argument. We establish that if an $L$-series with non-negative real coefficients (and positive coefficient at~$1$) converges in some half-plane, and has analytic continuation as an entire function, then its analytic continuation must take positive real values on all of~$\RR$ (\docs{LSeries.positive_of_differentiable_of_eqOn}). Assuming that $L(\chi, 1) = 0$, we obtain that $L(\chi, s) \zeta(s)$ extends to an entire function (the zero of $L(\chi, {\cdot})$ cancels the pole of~$\zeta$) and is an $L$-series with nonnegative coefficients, so it must be positive
  in particular at $s = -2$. But $\zeta(-2) = 0$, leading to a contradiction. This establishes the general non-vanishing result, \docs{DirichletCharacter.LFunction_ne_zero_of_one_le_re}.

  To prove Dirichlet's Theorem, we consider the $L$-series associated to the
  function that sends $n$ to~$\Lambda(n)$ (where $\Lambda$ is the von Mangoldt function)
  when $n \equiv a \bmod q$ and to zero otherwise. We show in \docs{ArithmeticFunction.vonMangoldt.LSeries_residueClass_eq} that this $L$-series is a linear combination of logarithmic derivatives of $L$-series associated to
  the Dirichlet characters mod~$q$.
  This step in the proof uses the orthogonality relations
  for (Dirichlet) characters, which we also contributed to Mathlib.
  The non-vanishing result above then implies that these logarithmic derivatives
  extend continuously to $\Re(s) \ge 1$, except for a simple pole at $s = 1$
  when $a$ and~$q$ are coprime. The existence of this simple pole
  then fairly quickly leads to a proof that $\sum_{n \equiv a \bmod q} \Lambda(n)/n$
  diverges, which in turn easily implies the existence of infinitely many
  prime numbers $p \equiv a \bmod q$. We state this final result in two versions,
  \docs{Nat.setOf_prime_and_eq_mod_infinite} (the set of these primes is infinite)
  and \docs{Nat.forall_exists_prime_gt_and_eq_mod} (for every natural number~$n$,
  there exists a larger prime $p \equiv a \bmod q$).

  The \href{https://leanprover-community.github.io/mathlib4_docs/Mathlib/NumberTheory/LSeries/Positivity.html}{positivity result} needs about 100~lines,
  \href{https://leanprover-community.github.io/mathlib4_docs/Mathlib/NumberTheory/DirichletCharacter/Orthogonality.html}{orthogonality} about 60~lines,
  the \href{https://leanprover-community.github.io/mathlib4_docs/Mathlib/NumberTheory/LSeries/Nonvanishing.html}{non-vanishing result} about 400~lines, and the
  \href{https://leanprover-community.github.io/mathlib4_docs/Mathlib/NumberTheory/LSeries/PrimesInAP.html}{final result} another 400~lines of Lean code.

\section{Acknowledgements}

 We would like to thank numerous members of the Lean community who assisted in the formalization project described here, in particular Chris Birkbeck, Riccardo Brasca, Kevin Buzzard, Johan Commelin and S\'ebastien Gou\"ezel. We also thank the anonymous referee for their close reading of the paper and perceptive comments.

 A key part of this project (the non-vanishing of $L$-functions on the abcissa of convergence) was written during a visit of Birkbeck and the second author to ETH Zürich in October 2024 as guests of the \emph{Forschunginstitut für Mathematik} (FIM), and we thank the FIM and its staff for their hospitality.

\printbibliography
% \begin{thebibliography}{dM+04}
% \bibitem[Buz20]{div0}
%  K. Buzzard, \emph{Division by zero in type theory: a FAQ}, Xena Project blog, 2020,
%  \url{https://xenaproject.wordpress.com/2020/07/05/division-by-zero-in-type-theory-a-faq/}

% \bibitem[DS05]{DS05}
%  F. Diamond and J. Shurman, \emph{A First Course in Modular Forms}. Graduate Texts in Mathematics vol. 228, Springer, 2005.

% \bibitem[Eb19]{Eb19}
%  M. Eberl, \emph{Nine Chapters of Analytic Number Theory in Isabelle/HOL}. In:
%  J. Harrison, J. O'Leary and A. Tolmach,
%  \emph{10th International Conference on Interactive Theorem Proving (ITP 2019)},
%  Leibniz International Proceedings in Informatics (LIPIcs), vol. 141, pp 16:1--16:19,
%  DOI: \href{https://doi.org/10.4230/LIPIcs.ITP.2019.16}{10.4230/LIPIcs.ITP.2019.16}.

% \bibitem[GK20]{GK20}
%  B. Gomes and A. Kontorovich, \emph{Riemann Hypothesis in Lean (with or without Mathlib)}, REU project, Rutgers University, 2020, \url{https://github.com/AlexKontorovich/Lean-RH}

% \bibitem[PNT+]{PNT+}
%  A. Kontorovich et al., \emph{The Prime Number Theorem And...}, \url{https://alexkontorovich.github.io/PrimeNumberTheoremAnd/web/}

% \bibitem[dM+24]{dM+24}
%  M. de Medeiros et al, \emph{Verified Foundations for Differential Privacy}, arXiv, December 2024, \url{https://arxiv.org/abs/2412.01671}

% \bibitem[Ri59]{Ri59}
%  B.~Riemann, \emph{\"Uber die Anzahl der Primzahlen unter einer gegebenen Gr\"osse [On the number of prime numbers less than a given size]}, Monatsberichte Berliner Akad., November 1859. (English translation by D.R. Wilkins: \url{https://www.claymath.org/wp-content/uploads/2023/04/Wilkins-translation.pdf})

% \bibitem[SY24]{SY24}
%  S. Song and B. Yao, \emph{Prime Number Theorem with Remainder Term}, Archive of Formal Proofs,
%  May 2024, \url{https://isa-afp.org/entries/PNT_with_Remainder.html}

% \end{thebibliography}
\end{document}